\documentclass{amsart}
\usepackage{amsmath,amssymb}

\newtheorem{theorem}{Theorem}

\newtheorem{lemma}[theorem]{Lemma}
\newtheorem{corollary}[theorem]{Corollary}

\theoremstyle{definition}

\def\drint{\ensuremath{- \hspace{-1.07em} \int}}
\def\rint{\ensuremath{- \hspace{-.9em} \int}}

\def\R{\mathbb R}
\def\Dom{\mathcal{D}}
\def\embed{\hookrightarrow}
\def\eps{\varepsilon}
\def\im{i}
\def\st{\,|\,}
\def\vt{\tau}
\DeclareMathOperator{\Ind}{ind}
\DeclareMathOperator{\Diff}{Diff}
\DeclareMathOperator{\spec}{spec}
\DeclareMathOperator{\sym}{\boldsymbol{\sigma}}

\begin{document}

\title{A note on the index of cone differential operators}
\author{Juan B. Gil}
\address{Department of Mathematics\\ Temple University\\
 Philadelphia, PA 19122}
\email{gil@math.temple.edu}
\author{Paul A. Loya}
\address{MIT, Department of Mathematics, Cambridge, MA 02139}
\email{ploya@alum.mit.edu}
\author{Gerardo A. Mendoza}
\address{Department of Mathematics\\ Temple University\\
 Philadelphia, PA 19122}
\email{gmendoza@math.temple.edu}


\begin{abstract}
We prove that the index formula for $b$-elliptic cone differential
operators given by Lesch in \cite{Le97} holds verbatim for operators
whose coefficients are not necessarily independent of the normal
variable near the boundary. We also show that, for index purposes,  the
operators can always be  considered on weighted Sobolev spaces.
\end{abstract}

\maketitle

\section{Introduction}
Let $M$ be a smooth compact manifold with
boundary, $\mathfrak m$ a smooth positive $b$-measure on $M$. With
respect to a suitable choice of a collar neighborhood $\pi:U\to
\partial M$ of the boundary and globally defined defining function
$x$ we may assume $\mathfrak m =\frac 1x dx \otimes \mathfrak
m_{\partial M}$ where $\mathfrak m_{\partial M}$ is a smooth positive
density on $\partial M$. Let $X$ be a vector field defined near
$\partial M$ such that $X$ is vertical with respect to $\pi$ and $Xx=1$.
Let $E$ be a vector bundle over $M$ and $A\in x^{-\nu}\Diff^m_b(M;E)$ be
$b$-elliptic, $\nu >0$ (see Lesch \cite{Le97} or Melrose \cite{RBM2} for
the notation). Fix a hermitian metric on $E$ and compatible connection
$\nabla$. Write $D_x=-\im\nabla_X$. The operator $A=x^{-\nu}P$ is said
to have coefficients independent of $x$ near the boundary if 
$xD_x P = PxD_x$ near $\partial M$.

Regard $A$ as an unbounded operator
 \begin{equation*}
 A:C_c^\infty(M;E)\subset x^{\mu}L^2_b(M;E)\to x^{\mu}L^2_b(M;E)
 \end{equation*}
 and denote by $\Dom_{\min}$ the domain of the
closure of $A$. It is convenient to assume $\mu=-\nu/2$; we can
always reduce to this case by conjugation with
$x^{\mu+\nu/2}$. In Section 2.4 of \cite{Le97} Lesch gives an analytic 
formula for the index of $A$ on $\Dom_{\min}$, assuming
that it has coefficients independent of $x$ near $\partial M$. 
This index formula is derived via heat trace asymptotics, also obtained 
in \cite{Le97}. In this note we show that the formula remains true also
without the assumption on the coefficients.

The results presented here represent a significant simplification of
various aspects of the analysis of $b$-elliptic cone operators and show
that simplifying assumptions made by various authors can be used in the
general case. For instance, operators with coefficients independent of
$x$ have a certain scalability property which was crucial in the works
of Cheeger \cite{Ch83} and Lesch \cite{Le97}.
This independence was also assumed by Schulze, Shatalov
and Sternin in \cite{SSS97} although, in the case they considered, a
much simpler homotopy invariance argument can be used to avoid this
assumption.

 From the technical point of view, the simplification comes about by
observing first that for the purpose of index calculations, one can
reduce to the case where the operator has coefficients independent of
$x$ near $\partial M$ (Theorem \ref{RedToConst}) and second, that one
can assume that the domain is the weighted Sobolev space
$x^{\nu/2-\eps}H^m_b(M;E)$ with suitable small $\eps>0$ rather than
$\Dom_{\min}$ (Theorem \ref{RedToSobolev}) by replacing $A$ with
$x^{\eps}A$. 

\section{The index formula}

Let $A\in x^{-\nu}\Diff_b^m(M;E)$ be $b$-elliptic. We shall need the
following lemma which also establishes the notation.

\begin{lemma} \label{NormEquiv}
On $\Dom_{\min}(A)$, with small enough $\eps>0$, the operator norm
\begin{equation*}
\|u\|_A = \|u\|_{x^{-\nu/2}L^2_b} + \|Au\|_{x^{-\nu/2}L^2_b}
\end{equation*} and the norm
\begin{equation*}
\|u\|_{A,\eps} = \|u\|_{x^{\nu/2-\eps}L^2_b} + \|Au\|_{x^{-\nu/2}L^2_b}.
\end{equation*} are equivalent.
\end{lemma}
\begin{proof} 
Recall that the embedding $x^{\nu/2-\eps}L^2_b\embed
x^{-\nu/2}L^2_b$ is continuous for $\eps<\nu$. The equivalence of the
norms follows from the continuity of
$(\Dom_{\min}(A),\|\cdot\|_A)\embed x^{\nu/2-\eps}L^2_b$ which is a
consequence of the closed graph theorem.
\end{proof}

Write $A=A_0+xA_1$ with $A_0$ having coefficients independent of $x$ 
near $\partial M$. Let $\omega\in C_c^\infty(\R)$, $\omega=1$ near $0$.
Furthermore, for $\vt>0$ let $\omega_\vt=\omega(x/\vt)$ and let
\begin{equation*}
 A^\vt=\omega_\vt A_0 + (1-\omega_\vt)A.  
\end{equation*}
Clearly, $A$ and $A^\vt$ have the same conormal symbol.

\begin{theorem}\label{RedToConst} 
For small enough $\vt$, $A^\vt$ is also $b$-elliptic, and therefore
$\Dom_{\min}(A^\vt)=\Dom_{\min}(A)$. Moreover, $A^\vt\to
A$ in the graph norm of $A$ as $\vt\to 0$. Thus on $\Dom_{\min}(A)$
\begin{equation*}
 \Ind A^\vt=\Ind A \quad\text{for every small } \vt. 
\end{equation*}
\end{theorem}
\begin{proof} 
Let $\sym(A)$ denote the totally characteristic principal
symbol of  $A$ order $m$. Then
\begin{align*}
\sym(A^\vt)
 &=\omega_\vt\sym(A_0) + (1-\omega_\vt)\sym(A)\\
 &=\omega_\vt\sym(A) + (1-\omega_\vt)\sym(A)- x\omega_\vt\sym(A_1)\\
 &=\sym(A)-\vt\tilde\omega_\vt\sym(A_1)
\end{align*} 
with  $\tilde\omega_\vt=(x/\vt)\omega(x/\vt)$. Since
$\tilde\omega_\vt$ is bounded, $\vt \tilde\omega_\vt$ is small for
$\vt$ small, and thus the invertibility of $\sym(A)$ implies that of
$\sym(A)-\vt\tilde\omega_\vt\sym(A_1)$  for such $\vt$. Hence
$A^\vt$ is $b$-elliptic too. Since $A$ and $A^\vt$ have the same
conormal symbol, part 1 of \cite[Proposition 4.1]{GiMe01} gives that
$\Dom_{\min}(A^\vt)=\Dom_{\min}(A)$. 

 From the $b$-ellipticity of $A$ it follows that there is a bounded
parametrix $Q:x^\gamma H^s_b\to x^{\gamma+\nu}H^{s+m}_b$ such that
\begin{equation*}
 R=I-QA:x^\gamma H^s_b\to x^{\gamma}H^\infty_b
\end{equation*} 
is bounded for all $s$ and $\gamma$. Write
\begin{align*}
 A-A^\vt = x\omega_\vt A_1
 &= x\omega_\vt A_1 QA + x\omega_\vt A_1 R \\
 &= \vt \tilde\omega_\vt A_1 QA + x\omega_\vt A_1 R.
\end{align*} 
Now, $A_1 Q: x^{-\nu/2}L^2_b\to x^{-\nu/2}L^2_b$ is bounded, so if
$u\in \Dom_{\min}(A)$, then
\begin{equation*}
 \|\vt \tilde\omega_\vt A_1 QAu\|_{x^{-\nu/2}L^2_b} \leq
   c\,\vt\|Au\|_{x^{-\nu/2}L^2_b} \leq c\,\vt\|u\|_A.
\end{equation*} 
Write $x\omega_\vt A_1 R= \vt^{1-\eps}(\frac x{\vt})^{1-\eps}
\omega_\vt\,x^{\eps}A_1 R$ and note that
\begin{equation*}
 x^{\eps}A_1 R: x^{\nu/2-\eps}L^2_b\to x^{-\nu/2}L^2_b
\end{equation*} 
in continuous. Then, using Lemma~\ref{NormEquiv} we get
\begin{equation*}
 \|x \omega_\vt A_1 Ru\|_{x^{-\nu/2}L^2_b}
 \leq \tilde c\,\vt^{1-\eps}\|u\|_{x^{\nu/2-\eps}L^2_b}
 \leq c\,\vt^{1-\eps}\|u\|_A.
\end{equation*} 
Altogether,
\begin{equation*}
 \|(A-A^\vt)u\|_{x^{-\nu/2}L^2_b} \leq C\,\vt^{1-\eps}\|u\|_A
\end{equation*} 
and thus $A^\vt\to A$ as $\vt\to 0$.
\end{proof}
                      
Using Lesch's formula \cite[Cor.~2.4.7]{Le97} and Proposition~3.14 in
\cite{GiMe01} we get
\begin{corollary} \label{cor:index} 
If $A_\Dom$ is an arbitrary closed extension of $A$, then
\begin{equation*}
\Ind A_\Dom=\drint_{M} \omega_A + \frac{1}{\nu}
\hat{\eta}(A)+\dim\Dom/\Dom_{\min}.
\end{equation*}
\end{corollary} 
The various terms occurring in the index formula are defined as follows. 
First, $\omega_A$ denotes the local index density of $A$. One way to 
define this density is as the constant term in the
difference of the small time fiberwise trace asymptotics of the heat
operators for $A^*_\Dom A_\Dom$ and $A_\Dom A^*_\Dom$. Actually,
$\omega_A$ is determined from the homogeneous terms in
$A$, and hence is defined independent of $\Dom$. In general,
$\omega_A$ diverges like $1/x$ at $\partial M$. The integral
$\rint_M \omega_A$ denotes the regularized integral of $\omega_A$, 
cf.\ \cite[Def.\ 2.1.3]{Le97}, and is defined as the usual integral of
$\omega_A$ off the collar $\pi:U\to \partial M$ of $\partial M$; and 
on the collar, $\rint_U \omega_A$ is defined as $\int_U \omega_1$, 
where $\omega_A = \pi^*\omega_0 \frac 1x dx + \omega_1$ on $U$, with 
$\omega_0$ a smooth density on $\partial M$ and with $\omega_1$ a 
smooth density on $U$. If $A$ has coefficients independent of $x$ near
$\partial M$, then $\omega_A$ vanishes identically near $\partial M$. 
In particular, in this case, $\rint_M \omega_A = \int_M \omega_A$ is just 
the usual integral of $\omega_A$.

In order to define $\hat{\eta}(A)$ let us first recall the definition of
the indicial operator of $A$. On the collar, we can write 
$A = x^{-\nu} \sum_{k=0}^m P_{m-k}(x) (xD_x)^k$, where $P_{m-k}(x)$ is
a  family of differential operators on $\partial M$ depending smoothly
on $x$. The indicial operator of $A$ is defined as the operator
\begin{equation*}
 A_\wedge = x^{-\nu}\sum_{k=0}^m P_{m-k}(0) (xD_x)^k 
\end{equation*}
on the model cone $N^\wedge= \overline{\R}_+ \times \partial M$. The eta
invariant $\hat{\eta}(A)$ is defined as the regular value at $z=0$ of
\begin{equation}
\Gamma(z)\big( \hat\zeta({A^*_{\wedge,\min}A_{\wedge,\min}},z)
-\hat\zeta({A_{\wedge,\min}A^*_{\wedge,\min}},z) \big),
\label{eta}
\end{equation} 
where $A_{\wedge,\min}$ denotes the closure of $A_\wedge$,
$\Gamma(z)$ is the gamma function, and where the zeta functions are
defined in the following way. Take for instance
$L=A^*_{\wedge,\min} A_{\wedge,\min}$. In Section 2.2 of \cite{Le97},
it  is shown that
\begin{equation*}
 k= \int_{\partial M}\mathrm{tr}(e^{-tL}(x,p,x,p))\mathfrak m_{\partial M}(p),
\end{equation*} 
is a function only of the singular coordinate $s = t/x^\nu$. Then,
$\hat\zeta(L,z)$ is defined as the transform
\begin{equation*}
 \hat\zeta(L,z):= \frac{1}{\Gamma(z)} \int_0^\infty s^{z-1}k(s)ds.
\end{equation*} 
By Proposition~2.2.6 of \cite{Le97}, $\hat\zeta(L,z)$ is meromorphic on
a half plane $\Re z < \delta$ for some $\delta > 0$, with only possible
simple poles at the points $(\dim M -k)/m$, $k=0,1, \ldots$. Note that
\`{a} priori, the function defined by \eqref{eta} may have a simple
pole at $z=0$; however, as shown in \cite[Cor.\ 2.4.7]{Le97}, this function 
is in fact regular there.

As already mentioned, the index formula of Corollary~\ref{cor:index} 
is an extension of Lesch's result \cite[Cor.~2.4.7]{Le97} to operators 
which do not necessarily have coefficients independent of $x$ near the
boundary. This index formula also generalizes the theorem of
Br\"{u}ning and Seeley \cite{BrSe88} for first order regular singular
operators. Index formulas for cone differential operators were first
proved by Cheeger \cite{Ch83} for the Gauss-Bonnet and signature
operators on a conic manifold, and were later generalized by Chou
\cite{Cho85} to  Dirac operators. For formulas in the
pseudodifferential situation see Fedosov, Schulze and Tarkhanov
\cite{FST99}. For formulas in the totally characteristic case ($\nu=0$), 
see Melrose \cite{RBM2} and Piazza \cite{Pi93}. We should point out, 
however, that the results of this note only apply when $\nu>0$. 

\section{Index on Sobolev spaces}

In general, $\Dom_{\min}$ is not a Sobolev space. The problem lies in
the possible presence of elements of $\spec_b(A)$ along the line
$\Im\sigma = -\nu/2$. However, for index purposes, one can conveniently
reduce the analysis to a slightly modified operator whose closure has a
Sobolev space as its domain:

\begin{theorem}\label{RedToSobolev} 
Let $A$ be $b$-elliptic. Let $A_\eps=x^{\eps}A$, and regard it as an 
unbounded operator on $x^{-(\nu-\eps)/2}L^2_b(M;E)$.
If $\eps>0$ is sufficiently small, then
\begin{equation*}
 A_\eps: x^{(\nu-\eps)/2}H^m_b(M;E)\to x^{-(\nu-\eps)/2}L^2_b(M;E)
\end{equation*}
is Fredholm, and
\begin{equation*}
 \Ind A_\eps = \Ind(A,\Dom_{\min}(A)).
\end{equation*}
\end{theorem}
\begin{proof} 
Write $A=x^{-\nu}P$ with $P\in \Diff^m_b(M;E)$. Let $\eta>0$ be so small 
that there is no $\sigma\in \spec_b(A)$ with $\nu/2-\eta\leq
\Im\sigma<\nu/2$ or $-\nu/2<\Im\sigma\leq  -\nu/2+\eta$. The kernel of
$A$ on tempered distributions $x^{-\infty}H^{-\infty}_b(M;E)$ is the
same as that of $P$, which we'll denote $K(P)$. Recall that
$\Dom_{\max}(A)=\{u\in x^{-\nu/2}L^2_b\st A u\in x^{-\nu/2}L^2_b\}$. 
The kernel $K_{\max}(A)$ of $A:\Dom_{\max}(A)\subset x^{-\nu/2}L^m_b\to
x^{-\nu/2}L^m_b$ consists those elements of $K(P)$ whose Mellin transforms 
are holomorphic in $\Im\sigma\geq \nu/2$ since these elements of $K(P)$ 
belong to $x^{-\nu/2}L^2_b$ and $Au\in x^{-\nu/2}L^2_b$ trivially. That is,
their Mellin transforms are holomorphic on $\Im\sigma>\nu/2-\eta$. Thus
$K_{\max}(A)=K_{\max}(A_\eps)$ if $0<\eps<\eta$.  On the other hand,
the kernel $K_{\min}(A)$ of $A:\Dom_{\min}(A)\subset x^{-\nu/2}L^m_b\to
x^{-\nu/2}L^m_b$ consists those elements of $K(P)$ whose Mellin transforms 
are holomorphic in $\Im\sigma> -\nu/2$; indeed in part 1 of 
\cite[Proposition 3.6]{GiMe01} it is shown show that
$\Dom_{\min}(A)=\Dom_{\max}(A)\cap x^{\nu/2-\eta}H^m_b$. Thus if
$\eps<\eta$ then $K_{\min}(A)=K_{\min}(A_\eps)$. 
Thus $\dim K_{\min}(A)=\dim K_{\min}(A_\eps)$. 

Finally, note that the formal adjoint of $A$
in $x^{-\nu/2}L^2_b$ is $A^\star = x^{-\nu}P^\star$, where $P^\star$ is
the formal adjoint of $P$ in $L^2_b$, and likewise
$A_\eps^\star = x^{-\nu+\eps}P^\star$. Now recall that the Hilbert
adjoint of $A|_{\Dom_{\min}}$ is $A^\star$ with domain
$\Dom_{\max}(A^\star)$ so the first part of the argument yields
$\dim K_{\max}(A^\star)=\dim K_{\max}(A_\eps^\star)$.
\end{proof}

\end{document}